\documentclass[11pt,a4paper]{article}
\usepackage{amsfonts}                 
\usepackage{amssymb}
\usepackage[latin1]{inputenc}         
\usepackage{color}
\usepackage{tikz}
\usetikzlibrary{arrows, decorations.markings, decorations.pathmorphing, backgrounds, positioning, fit, petri}

\newcommand{\prae}[2]{\mbox{$\langle #1\mid #2\rangle $}}

\newcommand{\qed}{\hspace*{\fill}$\Box $\\ \vspace{0.3cm}} 
\newcommand{\bewanf}{\noindent {\bf Proof:}\quad }    
\newcommand{\map}[2]{$f\colon #1\to #2$}
\newcommand{\satzanf}{\begin{samepage}\begin{satz}}
\newcommand{\satzende}{\end{satz}\end{samepage}}
\newcommand{\koranf}{\begin{samepage}\begin{kor}}
\newcommand{\korende}{\end{kor}\end{samepage}}
\newcommand{\lemanf}{\begin{samepage}\begin{lem}}
\newcommand{\lemende}{\end{lem}\end{samepage}}
\newcommand{\bspanf}{\begin{beispiel}}
\newcommand{\bspende}{\end{beispiel}}
\newcommand{\defanf}{\begin{definition}}
\newcommand{\defende}{\end{definition}}
\newcommand{\propanf}{\begin{proposition}}
\newcommand{\propende}{\end{proposition}}
\newtheorem{satz}{Theorem}[section]             
\newtheorem{lem}[satz]{Lemma}
\newtheorem{beispiel}[satz]{Example}
\newtheorem{definition}[satz]{Definition}
\newtheorem{proposition}[satz]{Proposition}
\newtheorem{kor}[satz]{Corollary}
\newcommand{\be}{\begin{equation}}
\newcommand{\ee}{\end{equation}}
\newcommand{\bl}[1]{\begin{equation}\label{#1}}
\newcommand{\ben}{\begin{enumerate}}
\newcommand{\een}{\end{enumerate}}
\begin{document}
\title{Combinatorial Relative Asphericity}
\author{Jens Harlander and Stephan Rosebrock}
\maketitle
\thispagestyle{empty}               
\begin{abstract} Relative notions of combinatorial asphericity have been used to prove that injective labeled oriented trees (which encode spines of ribbon 2-knots) are aspherical. This article presents an overview and comparison of the different notions of relative combinatorial asphericity. It also contains new results concerning characterizations of relative DR and tests that imply relative combinatorial asphericity. The last section of the article is devoted to examples that illustrate the concepts and the use of the tests given.
\end{abstract}


\section{Introduction}

In \cite{HR01} Huck and the second author showed that prime injective labeled oriented trees (which encode spines of ribbon 2-knots) are aspherical. We realized that in order to extend this result to all injective labeled oriented trees we needed notions of combinatorial asphericity relative to a sub-LOT. Relative vertex asphericity (VA) is used in \cite{HaRo17} to prove that injective labeled oriented trees are aspherical. A simpler proof that gives a stronger result is presented in \cite{HaRo19b}. 

The purpose of this article is to provide an overview and comparison of the different notions of combinatorial asphericity relative to a subcomplex. Given a pair of 2-complexes $K\subseteq L$, we say that $L$ is diagrammatically reducible (DR) relative to $K$ if every spherical diagram over $L$ can be reduced to a spherical diagram over $K$ by folds across edges with 2-cells from $L-K$. If folds across vertices are also allowed then this gives vertex asphericity (VA) relative to $K$. Precise definitions can be found in Section 2. Either notion implies that all second homotopy of $L$ is concentrated in $K$: The inclusion induced map $\pi_2(K)\to \pi_2(L)$ is surjective. If in addition one also wants $\pi_1$-injectivity $\pi_1(K)\to \pi_1(L)$ one needs a stronger version of relative combinatorial asphericity, a concept we called directed DR away from $K$. 

Bogley and Pride defined relative presentations $\langle H, {\bf x}\ |\ {\bf r} \rangle$, where $H$ is a group and ${\bf r}$ is a set of words in $F({\bf x})*H$.  See \cite{BP92} and \cite{BEW18}. They also define diagrammatic reducibility for relative presentations. Relative presentations have a more group theoretic flavor, but it turns out that DR for relative presentations is closely related to our concept of directed DR.\\


The paper is organized as follows. In Section 2 we present the different notions of relative combinatorial asphericity. In Section 3 we explain how the various notions relate to each other. In Section 4 we prove a characterization of relative diagrammatic reducibility along the lines of Corson-Trace. In Section 5 we present weight tests and max/min results with which one can show relative combinatorial asphericity. Results that are already published elsewhere are presented without proofs. This section contains new results that have not been published. The last Section 6 is devoted to examples that illustrate the concepts and the use of the tests given in Section 5. Most examples are from the class of labeled oriented trees.

\section{Notions of Combinatorial Relative Asphericity}

A map $f\colon X\to Y$ between complexes is {\em combinatorial} if $f$ maps open cells of $X$ homeomorphically to open cells of $Y$. A {\em surface diagram} over a 2-complex $K$ is a combinatorial map \map{C}{K}, where $C$ is a surface with a cell structure. If $C$ is a 2-sphere we call $f$ a {\it spherical diagram}. Note that if we orient the cells in $C$ and label each cell $c$ of $C$ by $f(c)$, the labeling on $C$ carries all information of $f$. We refer to such a labeled 2-sphere also as a spherical diagram over $K$. If $K$ is non-aspherical, then there exists a spherical diagram which realizes a nontrivial element of $\pi_2(K)$. In fact, $\pi_2(K)$ is generated by spherical diagrams. So in order to check whether a 2-complex is aspherical or not it is enough to check spherical diagrams.\\

Let $K$ be a 2-complex. The link of a vertex $v$, lk$(K,v)$, is the boundary of a regular neighborhood of $v$ in $K$. So lk$(K,v)$ is a graph whose edges are the corners of 2-cells at $v$. Suppose $K$ is a standard 2-complex with a single vertex $v$ and oriented edge set $X$. Then the vertices of lk$(K,v)=$ lk$(K)$ are $\{ x^+, x^- \mid x\in X\}$, where $x^+$ is a point of the oriented edge $x$ close to the beginning, and $x^-$ is a point close to the ending of that edge.  The {\it positive link} lk$^+(K)$  is the {\it full subgraph} on the vertex set  $\{ x^+ \mid x\in X\}$
and the {\it negative link} lk$^-(K)$ is the full subgraph on the vertex set $\{x^- \mid x\in X\}$.

Let $f\colon S\to K$ be a surface diagram and $v\in S$ a vertex. Restricting to the link we obtain a combinatorial map  $f|_{\mbox{lk}(S,v)}\colon \mbox{lk}(S,v)\to\mbox{lk}(K)$ for every vertex $v\in S$ and we let $z(v)=c_1\ldots c_q$ be the image, which is a {\it cycle} (a closed edge path) in lk$(K)$. 

\pagebreak

\defanf Let $\Gamma$ be a graph and $\Gamma_0$ be a subgraph. Let $z=e_1...e_q$ be a cycle. We say
\ben
\item $z$ is {\em homology reducible} if it contains a pair or edges $e_i,e_j$ such that $e_i=\bar e_j$ (the bar indicates opposite orientation) and {\em homology reduced} otherwise. $z$ is called {\em reducible} if in addition $j=i+1,\enspace (\mbox{mod }q)$.
\item $z$ is {\em homology reducible relative to $\Gamma_0$} if there is pair of edges $e_i,e_j$ such that $e_i=\bar e_j$ which is contained in $\Gamma-\Gamma_0$ and {\em homology reduced relative to $\Gamma_0$} otherwise. $z$ is called  {\em reducible relative to $\Gamma_0$}  if in addition $j=i+1,\enspace (\mbox{mod }q)$.
\een
\defende

Let $f\colon C\to L$ be a spherical diagram. A vertex $v\in C$ is called a {\it folding vertex} if $z(v)=c_1\ldots c_q\in \mbox{lk}(L)$ is homology reducible. In that case the pair of 2-cells $(d_i, d_j)$ of $C$ containing the preimages of $c_i$ and $c_j$, respectively, satisfying $c_i=\bar c_j$ is called a {\it folding pair}. If $z(v)$ is reducible, then $j=i+1$ and $d_i$ and $d_j$ share a common edge in $C$ which is called a {\it folding edge}.

$f$ is called {\it vertex reduced} if it does not have a folding vertex and $f$ is called {\it reduced} if it does not have a folding edge. A 2-complex $L$ is called {\it vertex aspherical} (VA) if each spherical diagram over $L$ has a folding vertex. $L$ is called {\it diagrammatically reducible} (DR) if each spherical diagram over $L$ has a folding edge. The vertices in the boundary of a folding edge are folding vertices so DR implies VA. Certainly VA implies asphericity.

\defanf Let $K$ be a subcomplex of the 2-complex $L$.  We say that 
\begin{itemize}
\item {\em $L$ is VA relative to $K$} if every spherical diagram $f\colon C\to L$, $f(C)\not\subseteq K$, has a folding vertex with folding pair of 2-cells in $L-K$. 
\item {\em $L$ is DR relative to $K$} if every spherical diagram $f\colon C\to L$, $f(C)\not\subseteq K$, has a folding edge with folding pair of 2-cells in $L-K$. 
\end{itemize}
\defende

\satzanf\label{sVAsub} Let $L$ be a 2-complex and $K$ a subcomplex. If $K$ is VA and $L$ is VA relative to $K$ then $L$ is VA.\satzende

\noindent Proof. Assume $f\colon C\to L$ is a vertex reduced spherical diagram. Since $L$ is VA relative to $K$ we have that $f(C)\subseteq K$. So $f\colon C\to K$ is a vertex reduced spherical diagram, contradicting the assumption that $K$ is VA. \qed

\satzanf\label{sDRsub} Let $L$ be a 2-complex and $K$ a subcomplex. If $K$ is DR and $L$ is DR relative to $K$ then $L$ is DR.\satzende

\noindent Proof. Assume $f\colon C\to L$ is a reduced spherical diagram. Since $L$ is DR relative to $K$ we have that $f(C)\subseteq K$. So $f\colon C\to K$ is a reduced spherical diagram, contradicting the assumption that $K$ is DR. \qed

\pagebreak

It is clear that if $L$ is DR relative to $K$ then $L$ is VA relative to $K$.

\satzanf\label{sVA} If $L$ is VA relative to $K$, then $\pi_2(L)$ is generated, as $\pi_1(L)$-module, by the image of $\pi_2(K)$ under the map induced by inclusion. In particular, if $K$ is aspherical, then so is $L$.
\satzende

\noindent Proof. Every vertex reduced spherical diagram $f\colon C\to L$ has its image $f(C)$ in $K$. Thus $f$ represents an element in $\pi_2(K)$. Since $\pi_2(L)$ is generated by vertex reduced spherical diagrams, it follows that $\pi_2(L)$ is generated by the image of $\pi_2(K)$. \qed\\

There is another notion of relative combinatorial reducibility which is studied in \cite{HaRo19} by the authors. For a set $X$ call a subset $Y$ {\it proper} if $Y\ne X$ ($Y$ may be empty). Remember that if $K$ is a 2-complex with edge set $X$ and if $f\colon C\to K$ is a spherical diagram then an edge $e\in C$ is labeled by $x\in X$, if $f(e)=x$.

\defanf\label{dDirDR} Let $K$ be a 2-complex with edge set $X$. Let $Y$  be a  proper subset of $X$. We say that $K$ is
\begin{itemize}
\item {\em DR directed away from $Y$} if every spherical diagram $f\colon C\to K$ that contains an edge with label from $X-Y$  also contains a folding edge with label from $X-Y$;
\item {\em DR in all directions} if every spherical diagram $f\colon C\to K$ that contains an edge labeled $x\in X$ also contains a folding edge with label $x$. Note that this implies that $K$ is DR directed away from all proper $Y\subset X$.
\end{itemize}
\defende

If $Y=\emptyset$ then DR directed away from $Y$ simply means DR. In \cite{HaRo19} directed DR is defined via presentations. If the 2-complex $K$ of Definition \ref{dDirDR} is the standard 2-complex given by a finite presentation then Definition \ref{dDirDR} is the same as the one in \cite{HaRo19}
given here:
If $P=\langle X\mid R\,\rangle$ is a presentation and $Y\subset X$ is proper we also say that {\it $P$ is DR directed away from $Y$} meaning that the standard 2-complex $K(P)$ build from $P$ is DR directed away from the 1-cells corresponding to $Y$.\\

Let $K$ be a 2-complex with edge set $X$ and let $Y$ be a proper subset of $X$. We define $K_Y$ to be the subcomplex of $K$ with edge set $Y$ containing exactly those 2-cells of $K$ with all boundary 1-cells from $Y$. The following theorem is Theorem 2.2 of \cite{HaRo19}. We include a proof for the convenience of the reader:

\satzanf\label{RelPiS} Let $K$ be a 2-complex with edge set $X$. Suppose that $K$ is DR directed away from the proper subset $Y\subset X$.  Then 
\begin{enumerate}
\item $\pi_2(K)$ is generated (as a $\pi_1(K)$-module) by the image of the inclusion induced map $\pi_2(K_Y)\to \pi_2(K)$; furthermore
\item every disc diagram $g\colon D\to K$ with boundary labeled by a word in $Y$, that contains a label from $X-Y$, has a folding edge with label from $X-Y$. Consequently, the inclusion induced map $\pi_1(K_Y)\to \pi_1(K)$ is injective.
\end{enumerate}
\satzende

\noindent Proof. Suppose $f\colon C\to K$ is a reduced spherical diagram. If $f(C)$ is not contained in $K_Y$ then $C$ contains an edge $e$ so that $f(e)\not\in Y$. Since we assumed that $K$ is DR directed away from $Y$ it follows that $C$ contains a folding edge $e'$ so that $f(e')\not\in Y$, contradicting the assumption that $f\colon C\to K$ is reduced. Since $\pi_2(K)$ is generated (as a $\pi_1(K)$-module) by reduced spherical diagrams, the first statement follows.

Suppose $g\colon D\to K$ is a disc diagram as in statement (2). We double $D$ and construct a spherical diagram $g'\colon C=D_1\cup D_2 \to K$, where $D_1$ is mapped by $g$ and $D_2$ is mapped by $-g$ (an orientation reversion followed by $g$).  Note that $C$ contains an edge with label not in $Y$. Since $K$ is DR away from $Y$ this spherical diagram contains a folding edge with label not from $Y$. This folding edge can not occur on $\partial D_1=\partial D_2$. Thus $D_1$ or $D_2$ contain an interior folding edge with label not in $Y$ and hence so does $D$.

We next show $\pi_1$-injectivity. Suppose $w$ is a word in $Y^{\pm 1}$ that represents a non-trivial element of $\pi_1(K_Y)$ that maps to a trivial element in $\pi_1(K)$. Then there exists a reduced Van Kampen diagram $f\colon M\to K$ where the boundary of $M$ is labelled by $w$ and is mapped to $K_Y$. Here $M$ is a planar simply connected region with a cell structure. Note that $M$ is a tree with discs attached at some vertices. One of these discs, say $\bar D$, contains an edge $e$ such that $f(e)\not\in K_Y$, otherwise $f(M)\subseteq K_Y$ which would imply that $w=1$ in $\pi_1(K_Y)$. Thus ${\bar g}=f|_{\bar D}\colon {\bar D}\to K$ is a disc diagram as in statement (2) and hence contains a folding edge, contradicting the fact that we assumed $f\colon M\to K$ is reduced.\qed

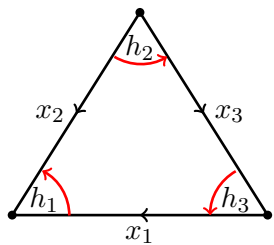
\begin{figure}[ht]\centering
\begin{tikzpicture}[scale=0.84]
 
\fill (0,0) circle (2pt);
\fill (4,0) circle (2pt);
\fill (2,3.2) circle (2pt);
 
\begin{scope}[decoration={markings, mark=at position 0.5 with {\arrow{>}}}]
\draw [postaction={decorate},line width=1pt] (4,0) -- (0,0) node[midway, below]{$x_1$};
\draw [postaction={decorate},line width=1pt] (2,3.2) -- (0,0) node[midway, left]{$x_2$};
\draw [postaction={decorate},line width=1pt] (2,3.2) -- (4,0) node[midway, right]{$x_3$};
\end{scope}
 
\node at (0.5,0.25) {$h_1$};
\node at (3.5,0.25) {$h_3$};
\node at (2,2.65) {$h_2$};
 
\draw [<-,red,line width=1pt ] (0.47,0.7) arc (62:0:0.8);
\draw [->,red,line width=1pt ] (1.6,2.5) arc (239:301:0.8);
\draw [<-,red,line width=1pt ] (3.1,0) arc (180:120:0.8);
 
\end{tikzpicture}
\caption{\label{arelrel} The relator $\bar r$ for  ${\hat r}=x_1h_1x_2^{-1}h_2x_3h_3$}
\end{figure}

We end this section with a brief mention of relative presentations due to Bogley and Pride \cite{BP92}. 
A {\em relative presentation} $\hat P=\langle H, {\bf x}\mid {\bf \hat r} \rangle$ consists of a group $H$, a generating set ${\bf x}$ and relator set ${\bf \hat r}\subseteq H*F({\bf x})$ ($F({\bf x})$ is the free group generated by $\bf x$). Bogley and Pride defined diagrammatic reducibility for relative presentations in terms of pictures rather than diagrams. Pictures and diagrams are dual concepts. Here is the idea in terms of marked diagrams.
Let $\bar P$ be the presentation obtained from $\hat P$ by forgetting all the $H$-information. That is $\bar P=\langle {\bf x}\mid  \bar {\bf r} \rangle$ where $\bar r\in \bar {\bf r}$ is the image of $\hat r\in {\bf \hat r}$ under the projection $H*F({\bf x})\to F({\bf x})$. The 2-complex $K(\hat P)$ is the 2-complex $K(\bar P)$ with corners marked by elements from $H$: If $\hat r=x_{1}h_1x_{2}h_2\ldots x_{t}h_t\in {\bf \hat r}$, where $x_{i}\in {\bf x}^{\pm 1}$ and $h_i\in H$, then $\bar r=x_{1}\ldots x_{t}$ and we mark the oriented corner from the $i$-th edge to the $(i+1)$-th edge in the boundary of the 2-cell $\Delta(\bar r)$ (oriented clockwise) by $h_i$. An example is given in Figure \ref{arelrel}.

In  a surface diagram $\hat f\colon \hat F \to K(\hat P)$ the corners in the 2-cells are marked by elements from $H$. The corner marking is obtained by pulling back the corner marking of the cells of $K(\hat P)$. We can assign group elements $h(v)$ to the vertices $v$ in the diagram in the following way: If $c_1(v)\ldots c_l(v)$ is the clockwise corner cycle (or path, in case $v$ is a boundary vertex) at the vertex $v$, and $c_i(v)$ is marked with $h_i$, then $h(v)=h_1\ldots h_l$. Note that $h(v)$ is defined only up to cyclic permutation in case $v$ is an interior vertex.

A spherical diagram $\hat f\colon \hat C \to K(\hat P)$ is {\em admissible} if $h(v)=1$ for all but possibly one vertex $v_0$ of $\hat C$. A relative presentation is defined to be {\it diagrammatically reducible DR} if every admissible spherical diagram has a folding edge. 

DR for relative presentations is important in the study of equations over groups. It is known that if $\hat P$ is DR, then the map $H\to G(\hat P)$ is injective, that is the set of equations $\bf \hat r=1$ over $H$ has a solution in an overgroup. For a good overview on the topic of relative presentation we refer to \cite{BEW18}.

\section{Comparing the notions of relative asphericity}

\propanf Let $K$ be a 2-complex with edge set $X$. Let $Y$  be a  proper subset of $X$. If $K$ is DR directed away from $Y$, then $K$ is DR relative to $K_Y$.
\propende

Proof. Assume $K$ is not DR relative to $K_Y$. Then there exists a spherical diagram $f\colon C\to K$ such that $f(C)\not\subseteq K_Y$ where all pairs of 2-cells which may be reduced lie in $K_Y$. Since $f(C)\not\subseteq K_Y$ we have that $C$ contains an edge labelled by an element of $X - Y$. Since  all pairs of 2-cells which may be reduced lie in $K_Y$ we have that all folding edges of $C$ are labelled by elements of $Y$. So $K$ is not DR directed away from $Y$.\qed

If $K$ is DR relative to $K_Y$ then $K$ does not have to be DR directed away from $Y$. This is because a pair of cancelling 2-cells in a spherical diagram might map to $K - K_Y$ but the common edge of these 2-cells maps to $Y$.

\bspanf {\em Let $P$=\prae{a,b,c}{bac^{-1},cb^{-1}a^{-1}} be a presentation. $K(P)$ is the torus. There is a disk diagram $D$ with boundary reading $aba^{-1}b^{-1}$ achieved by gluing the two relator disks along $c$ together. Glue $D$ to $-D$ to obtain a spherical diagram which is reducible at $a$ and $b$ only. This shows that $P$ is not DR away from $Y=\{ a,b\}$.

Since $K_Y$ is the 2-complex modeled on the presentation $\langle a,b\mid\ \rangle$ which has no relators we have that $K(P)$ is DR relative $K_Y$ if and only if $K(P)$ is DR. But this presentation of the torus is certainly DR, so $K(P)$ is DR relative $K_Y$.}
\bspende

Let $K$ be the standard 2-complex given by the presentation $P=\langle X \mid R\,\rangle$ and let $Y\subseteq X$ proper. We can associate to this a relative presentation: Let $H=\pi_1(K_Y)$,  ${\bf x}=X - Y$ and ${\bf \hat r}$ be the set obtained from $R$ in the following way: We have a homomorphism $\phi\colon F(X)\to H*F({\bf x})$ by sending $x\in Y$ to the group element in $H$ that it presents, and $x$ to $x$ if $x\in X - Y$. Let $\hat r=\phi(r)$, $r\in R$. Then we have a relative presentation $\hat P=\langle H, {\bf x}\mid {\bf \hat r} \rangle$.  The following result, obtained by the authors, can be found in \cite{HaRo19}.

\satzanf $P$ is DR directed away from $Y$ if and only if the relative presentation $\hat P=\langle H, {\bf x}\mid {\bf \hat r} \rangle$ is DR.
\satzende

\vspace*{0.4cm}

Let $K$ be a 2-complex with edge-set $X$ and $Y\subset X$ a proper subset. Let $H=\pi_1(K_Y)$ and ${\bf x}=X-Y$. The following diagram visualizes the connections between the different notions of relative combinatorial asphericity:\\

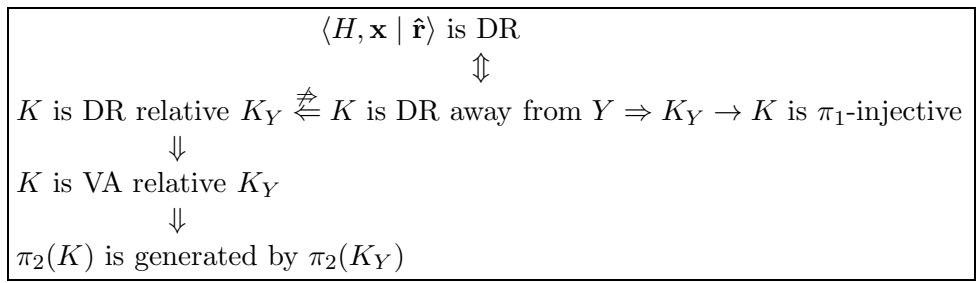
\fbox{\parbox{12.5cm}{\hspace*{4cm}$\langle H ,{\bf x} \mid {\bf \hat r} \rangle$ is DR\\
\hspace*{6cm}$\Updownarrow $\\
$K$ is DR relative $K_Y\stackrel{\not\Rightarrow}{\Leftarrow} K$ is DR away from $Y\Rightarrow K_Y\to K$ is $\pi_1$-injective\\
\hspace*{2cm}$\Downarrow $\\
$K$ is VA relative $K_Y$\\
\hspace*{2cm}$\Downarrow $\\
$\pi_2(K)$ is generated by $\pi_2(K_Y)$}}

\vspace{0.5cm}

\pagebreak

\section{A Characterization of Relative Diagrammatic\\ Reducibility}

For a complex $K$ let $\tilde K$ be its universal covering. Let $K^{(1)}$ be the 1-skeleton of $K$.
Corson and Trace \cite{CT00} have shown:

\satzanf\label{satz:CT}  The 2-complex $K$ is DR if and only if every finite subcomplex of the universial covering $\tilde K$ collapses into $\tilde K^{(1)}$.
\satzende

The authors have generalized this result to directed DR in \cite{HaRo19}. Here is the corresponding result for relative DR.

\satzanf\label{satz:CTgen} Let $L$ be a 2-complex and $K$ a subcomplex. Then $L$ is DR relative to $K$ if and only if every finite subcomplex of $\tilde L$ collapses into $p^{-1}(K)\cup \tilde L^{(1)}$, where $p\colon \tilde L\to L$ is the covering projection.
\satzende

Recall that an edge in a 2-complex is called {\em free} it it occurs exactly once in the boundary of exactly one 2-cell. A 2-complex is called {\it closed} if it does not have a free edge.

\begin{lem}\label{lem:CT} Let $K$ be a finite 2-complex and $d$ be a 2-cell in $K$. If $K$ is closed then there exists a reduced surface diagram $f\colon F\to K$ such that $F$ is closed and orientable and $d$ is contained in  $f(F)$.
\end{lem}

This lemma is in Corson-Trace \cite{CT96}, Theorem 2.1, stated without the fixed 2-cell $d$. The fact that $f$ hits a specified 2-cell will be important in the relative case. A detailed proof of Lemma \ref{lem:CT} can be found in \cite{HaRo19}, Lemma 3.2.\\

Proof (of Theorem \ref{satz:CTgen}). Assume first that $L$ is DR relative to $K$. Suppose the claim is false. Among all finite subcomplexes of $\tilde L$ that do not collapse into $p^{-1}(K)\cup \tilde L^{(1)}$ let $M$ be one with the minimal number of 2-cells. Note that $M$ does not have a free edge, because a collapse could be performed at that free edge to produce a complex with fewer 2-cells, contradicting minimality. 

Let $\tilde d$ be a 2-cell in $M$ not contained in $p^{-1}(K)\cup \tilde L^{(1)}$. It follows from Lemma \ref{lem:CT} that there exists a reduced surface diagram $\tilde f\colon  F\to M\subseteq \tilde L$, where $F$ is a closed orientable surface and $\tilde d$ is contained in $\tilde f(F)$. Let $f=p\circ \tilde f\colon F\to L$. Note that $d=p(\tilde d)$ is a 2-cell in $f(F)$ not contained in $K$.
We now proceed as in the proof of Lemma 2.1 in \cite{CT00}: Attach Van Kampen diagrams along cutting curves of $F$ to produce a simply connected 2-complex $N_0$ and combinatorial maps $F\stackrel{\alpha_0}\to N_0\stackrel{\beta_0}\to L$ such that $\beta_0 \circ \alpha_0=f$. The Van Kampen diagrams exist because $f$ lifts to $\tilde f$, so every closed curve in $F$ maps to a closed curve in $L$ that is homotopically trivial. Note furthermore that $N_0$ is the 2-skeleton of a cell decomposition of the 3-sphere $S^3$. Let $N$ be a 2-complex with the minimal number of 2-cells satisfying the following conditions:
\begin{enumerate}
\item $N$ is a simply connected 2-skeleton of a cell-decomposition of $S^3$;
\item There exist combinatorial maps $F\stackrel{\alpha}\to N\stackrel{\beta}\to L$ such that $\beta \circ \alpha=f$.
\end{enumerate}
Corson-Trace show that the attaching maps of the 3-cells of $S^3$ utilize all 2-cells of $N$ and each attaching map $\gamma\colon C\to N$ results in a reduced spherical diagram $\beta\circ \gamma\colon C\to L$. Choose a 3-cell so that $\beta\circ \gamma(C)$ contains $d$. Then $\beta\circ \gamma\colon C\to L$ is a reduced spherical diagram that is not a diagram over $K$, contradicting the assumption that $L$ is DR relative to $K$.

For the other direction assume that $f\colon C\to L$ is a spherical diagram which is not already a diagram over $K$. We can lift it to a spherical diagram $\tilde f\colon C\to \tilde L$. Now $\tilde f(C)=M$ is a finite subcomplex of $\tilde L$ not already contained in $p^{-1}(K)\cup \tilde L^{(1)}$. Since $M$ collapses into $p^{-1}(K)\cup \tilde L^{(1)}$, it has a free edge $\tilde e$ in the boundary of a 2-cell not contained in $p^{-1}(K)$. Any edge $e \in C$ so that $\tilde f(e)=\tilde e$ is a folding edge. Thus $f\colon C\to L$ is not reduced relative to $K$. \qed

The following corollary can be used as a tool to show that a given presentation defines an infinite group.  

\koranf Let $P$ be a presentation of a finite group and $T$ a subpresentation. Then $K(P)$ is DR relative to $K(T)$ if and only if $K(P)$ collapses into $K(T)$.
\korende

For completeness we finish this section with the directed DR version of the Corson-Trace result from \cite{HaRo19}. Comparing the two theorems shows again the subtle differences between relative DR and directed DR.

\satzanf Let $K$ be a 2-complex with edge set $X$. Let $Y$  be a  proper subset of $X$. Then $K$ is DR directed away from $Y$ if and only if every finite subcomplex of $\tilde K$ collapses into $p^{-1}(K_Y)\cup \tilde K^{(1)}$, where only edges that map to $x\in X - Y$ under the covering projection are used as collapsing edges.
\satzende

\section{Methods for Showing Relative Combinatorial\\ Asphericity}

We first give two results which use tree-like conditions in the link-graph to show relative combinatorial asphericity.

Let $\Gamma $ be a graph and $\hat\Gamma=\Gamma_1\cup \ldots \cup \Gamma_n$ be a union of disjoint subgraphs.  We write $\Gamma /\hat\Gamma$ for the graph obtained from $\Gamma$ by collapsing each $\Gamma_i$ to a vertex.

\defanf\label{dreltree} Let $\Gamma$ be a graph and $\hat \Gamma=\Gamma_1\cup \ldots \cup \Gamma_n$ be a disjoint union of subgraphs. We say $\Gamma$ is a {\em forest relative to $\hat\Gamma$} if $\Gamma /\hat\Gamma$ has no cycles. $\Gamma $ is called a {\em tree relative to $\hat\Gamma$} if in addition $\Gamma$ is connected.
\defende

A subcomplex $K$ of a 2-complex $L$ is called {\it full}, if for every 2-cell $d\in L$ where all boundary cells are in $K$ we have $d\in K$. The attaching maps of 2-cells of a 2-complex $K$ are said to have {\it exponent sum 0} if for a given orientation of the 1-cells of $K$ the attaching map of each 2-cell $d$ satisfies: when traveling along the boundary of $d$ in clockwise direction, one encounters the same number of positive as of negative edges.

 The following theorem is shown in \cite{HaRo19b}:

\satzanf\label{sforest} Let $L$ be a finite 2-complex with one vertex and $K=K_1\vee\ldots\vee K_n\subseteq L$. We assume the attaching maps of 2-cells in $L$ have exponent sum $0$, and the $K_i$ are full. If lk$^+(L)$ is a forest relative to $lk^+(K)$ or  lk$^-(L)$ is a forest relative to lk$^-(K)$ then $L$ is VA relative to $K$. Furthermore, the inclusion induced homomorphism $\pi_1(K_i)\to \pi_1(L)$ is injective for every $i=1,\ldots ,n$.
\satzende

The following theorem is proved in \cite{HaRo19}:

\satzanf\label{sminmax} Let $K$ be a finite 2-complex with one vertex and edge set $X$.  We assume the attaching maps of 2-cells in $K$ have exponent sum $0$. Assume that lk$^+(K)$ or lk$^-(K)$ is a forest. Then for each $x\in X$, $K$ is DR directed away from $\{ x\} $.
\satzende

If $\omega $ is a real valued function on the set of edges of a graph $\Gamma $ and $z=e_1\ldots e_p$ is an edge path in $\Gamma$, we write 
\[ \omega (z)=\sum_{1\le i\le p}\omega (e_i).\]

$\omega $ is called a {\it weight function}. The following is a weight test for directed diagrammatic reducibility. It is a generalized version of Gerstens weight test (see \cite{Ger87}). Its proof is given in \cite{HaRo19}.

\satzanf\label{sWT} Let $K$ be a 2-complex with cyclically reduced attaching maps of 2-cells and edge set $X\cup Y$ with $X=\{ x_1,\ldots ,x_n\}$ and $Y=\{ y_1,\ldots ,y_p\}$.
Suppose we can assign weights $\omega(e)\ge 0$ to the edges $e$ of lk$(K)$, such that: 
\begin{enumerate}
\item If $e$ connects $y_i^\epsilon$ with $y_j^\delta$, ($\epsilon ,\delta =\pm$) then $\omega(e)\ge 1$; 
\item If one of $e$'s boundary vertices is $y_i^+$ or $y_i^-$, then $\omega(e)\ge 1/2$;
\item If $z$ is a reduced cycle in  lk$(K)$, then $\omega(z)\ge 2$;
\item Let $d$ be a 2-cell of length $\kappa (d)$ from $K$, then $\sum_{c\in d}\omega(c)\le \kappa (d)-2$.
\end{enumerate}
Then $K$ is DR directed away from $Y$.
\satzende

There are at least two weight tests which measure relative vertex asphericity. The most general is the following:

Let $L$ be a 2-complex and $K=K_1\vee\ldots\vee K_n $ be a subcomplex of $L$. Let $W(L,K)$ be lk$(L)$ where lk$(K_i)$ is replaced by two vertices $k^+_i,k^-_i$ and exactly one edge $e_i$ connecting these two vertices for each $1\le i\le n$. So the edges of $W(L,K)$ are the corners of 2-cells of $L$ not in a $K_i$ and edges $e_1,\ldots ,e_n$. The vertices of $W(L,K)$ are the vertices of lk$(L)$ without the vertices of all lk$(K_i)$ but with two additional vertices $k_i^+,k_i^-$ for each $K_i$.

\satzanf\label{sWTa} Let $L$ be a 2-complex with one vertex and cyclically reduced attaching maps of 2-cells and let $K=K_1\vee\ldots\vee K_n $ be a subcomplex of $L$ such that no edge of $K$ represents the trivial element in $\pi_1(K)$. If there is a weight function $\omega $ on the set of edges of $W(L,K)$ which satisfies:
\ben
\item $\sum_i \omega(c_i)\le q-2$ if $c_1,\ldots ,c_q$ are the corners of a 2-cell of $L-K$,
\item if $z$ is a homology reduced cycle in $W(L,K)$ then $\omega(z)\ge 2$,
\item $\omega (e_i)=0$ for $1\le i\le n$,
\een
then $L$ is VA relative $K$.
\satzende

\bewanf Assume $f\colon C\to L$ is a vertex reduced spherical diagram such that $f(C)\not\subset K$. Pull back the weights of $W(L,K)$ to corners in 2-cells of $C$ which map to $L-K$. Replace each maximal region $d\in C$ which is mapped to a single $K_i$ by a 2-cell $d'$ (we call $d'$ a {\it replaced 2-cell}) and achieve a new cell decomposition $C'$ of the 2-sphere. If $d$ is homeomorphic to a disk with $m$ holes then there are $m$ arcs in $C'$ bounding $d'$ on both sides.  Assign weight 0 to all  corners of $d'$. Assign the weights of $C$ to the corresponding corners of non-replaced 2-cells of $C'$.

The curvature of those 2-cells of $C'$ coming from 2-cells of $C$ mapped to $L-K$ have curvature less or equal to 0 by condition 1. A replaced 2-cell $d'\in C'$ has at least two corners by the condition that 
no edge of $K$ represents the trivial element in $\pi_1(K)$. It has weight $0$ and condition 1 is satisfied for this 2-cell, leading also to non-positive curvature for replaced 2-cells.

Assume $z'$ is the link of a vertex in $C'$. If $z'$ contains no corners of replaced 2-cells then it has weight at least two by condition 2 since $f$ is homology reduced. If $z'$ contains a corner of a replaced 2-cell, this corner will contribute 0 to the weight of $z'$. Since it does not appear in $W(L,K)$ if it is a corner in lk$^+(K_i)$ or lk$^-(K_i)$ or it has weight $0$ in $W(L,K)$ we have weight at least two by condition 2 for $z'$.

So we have non-positive curvature at vertices of $C'$ contradicting the Euler-Characteristic of $C'$ by the combinatorial Gauss-Bonet theorem.
\qed\\

Next we present a weight test with stronger requirements than the ones in Theorem \ref{sWTa} but on the other hand fits to more presentations satisfying the requirements. 
Assume $L$ is a 2-complex and $K=K_1\vee\ldots\vee K_n\subseteq L$. We assume $L$ contains a single vertex $v$. We define lk$(L,K)$, the link of $v$ in $L$ relative to $K$ in the following way: If $y_1,\ldots ,y_l$ are the edges of $K_i$ then we denote by $\Delta(K_i)$ the full graph on the vertices $y_i^{\pm 1}$ of lk$(K_i)$ together with an edge attached at each $y_i^+$ (a loop at that vertex) and at each $y_i^-$. Every pair of vertices in $\Delta(K_i)$ is connected by an edge, and at every vertex we have a loop. For each $i$ we remove lk$(K_i)$ from lk$(L)$ and insert $\Delta(K_i)$ instead. The resulting graph is lk$(L,K)$. The proof of the following theorem is given in \cite{HaRo19b}.

\satzanf\label{sWTb} Assume $K=K_1\vee\ldots\vee K_n\subseteq L$, each $K_i$ is full and the attaching maps of 2-cells of $K$ have exponent sum zero. Assume there is a weight function on the edges of lk$(L,K)$ satisfying
\begin{enumerate}
\item $\sum_i \omega(c_i)\le q-2$ if $c_1,\ldots ,c_q$ are the corners of a 2-cell of $L$ not contained in $K$,
\item if $z$ is a homology reduced cycle in $lk(L, K)$ containing at least one corner from  lk$(L,K)-\Delta(K)$ then $\omega(z)\ge 2$,
\item $\omega(c)=0$ if $c$ is an edge of $\Delta^{+}(K_i)$ or $\Delta^{-}(K_i)$ and\\
    $\omega(c)=1$ if $c$ connects a vertex of $\Delta^{+}(K_i)$ with one of $\Delta^{-}(K_i)$.
\end{enumerate}
Then $L$ is VA relative to $K$. If in addition the attaching maps of the 2-cells of $L$ have exponent sum zero, then all the inclusion induced homomorphisms $\pi_1(K_i)\to \pi_1(L)$ are injective.
\satzende

\vspace{1ex}

We recall the Freiheitssatz for 1-relator groups: Suppose $P=\langle x_1,\ldots ,x_n\mid r \rangle$ is a 1-relator presentation, where $r$ is a cyclically reduced word that contains all the generators. Then any proper subset $Y$ of $\{x_1,\ldots ,x_n\}$ generates a free subgroup of the corresponding group $G(P)$ with basis $Y$. The following three results, already contained in  \cite{HaRo19}, can be viewed as multi-relator versions of this celebrated result. For convenience to the reader we include short proofs.\\

If $P=\langle X\mid R\,\rangle$ is a presentation and $Y\subset X$ define $P_Y$ to be the subpresentation of $P$ with generators $Y$ and all relators that contain only generators of $Y$.

\satzanf Let $P=\langle x_1,\ldots ,x_n\mid r_1,\ldots ,r_m\rangle$ be a presentation such that $P$ is DR in all directions. Then the inclusion induced homomorphism $G(P_Y)\to G(P)$ is injective for every subset $Y$ of the generators. 
\satzende

\noindent Proof. If $Y$ is the set of generators of $P$ then $G(P_Y)=G(P)$ and the statement is true. If $Y$ is a proper subset of the set of generators then the result follows from Theorem \ref{RelPiS} using the fact that $P$ is DR directed away from $Y$.\qed

\pagebreak

\satzanf Let $P=\langle x_1,\ldots ,x_n\mid r_1,\ldots ,r_m\rangle$ and $Y$ a proper subset of the generators. Assume that each $r_i$ contains a generator not from $Y$. If $P$ is DR
directed away from $Y$, then $Y$ generates a free subgroup of $G(P)$ with basis $Y$.
\satzende

\noindent Proof. Since each $r_i$ contains a generator not from $Y$ we have that $P_Y=\langle Y\mid\ \rangle $ and $G(P_Y)$ is free. Now Theorem \ref{RelPiS} (2) gives the desired result.\qed

\koranf\label{cor:freiheit} Let $P=\langle x_1,\ldots ,x_n\mid r_1,\ldots ,r_m\rangle$ be a presentation where each $r_i$ contains all the generators. If $P$ is DR in all directions, then any proper subset $Y$ of $\{x_1,\ldots ,x_n\}$ generates a free subgroup of $G(P)$ with basis $Y$. 
\korende

Here is a strengthening of the classical Freiheitssatz which is proved in \cite{HaRo19}:

\satzanf\label{thm:1rel}Let $P=\langle x_1,\ldots ,x_n \mid r \rangle$ be a one-relator presentation of a group $G$ where $r$ is a cyclically reduced relator that is not a proper power. Then $P$ is DR in all directions.\satzende

\section{Examples}

Many examples we present are LOGs. A standard reference for labeled oriented graphs, LOGs for short, is \cite{Ro18}. Here are the basic definitions. A LOG is an oriented finite graph $\Gamma$ on vertices $\textbf{x}$ and edges ${\bf e}$, where each oriented edge is labeled by a vertex. Associated with a LOG $\Gamma $ is the {\em LOG-complex} $K(\Gamma)$, a 2-complex with a single vertex, edges in correspondence with the vertices of $\Gamma$ and 2-cells in correspondence with the edges of $\Gamma$. The attaching map of a 2-cell $d_e$ is the word  $xz(zy)^{-1}$, where $e$ is an edge of $\Gamma$ starting at $x$, ending at $y$, and labeled with $z$.

A labelled oriented graph is called {\em compressed} if no edge is labelled with one of its vertices. A LOG $\Gamma $ is called {\em boundary reducible} if there is a boundary vertex $x\in\Gamma$ which does not occur as an edge label and {\em boundary reduced} otherwise.
A LOG is {\em injective} if each vertex occurs as an edge label at most once. 
An injective LOG is called {\em reduced} if it is compressed and boundary reduced.
A {\em labeled oriented tree}, LOT, is a labeled oriented graph where the underlying graph is a tree. If $\Gamma$ is a LOT and $\Gamma_1$ is a sub-tree of $\Gamma$ with at least one edge, such that each edge label of $\Gamma_1$ is a vertex of $\Gamma_1$, then we call $\Gamma_1$ a {\em sub-LOT} of $\Gamma$.\\

A main result of \cite{HaRo19b} is:

\satzanf\label{sinj} Let $\Gamma$ be a compressed injective LOT. Then $K(\Gamma)$ is VA.\satzende

The main step used in the proof is the following: 
If $\Gamma $ is a reduced injective LOT which contains disjoint sub-LOTs $\Gamma_1,\ldots ,\Gamma_n$, then in most cases $K(\Gamma)$ is VA relative to $K(\Gamma_1)\cup\ldots\cup K(\Gamma_n)$. This is shown by using the weight test of Theorem \ref{sWTb}. Since by induction the $\Gamma_i$ are VA the result is shown by Theorem \ref{sVAsub}.\\

The next lemma gives a tool for proving directed DR in many situations (see \cite{HaRo19}).
Let $T_1=\langle x_1,\ldots ,x_k\mid r_1,\ldots ,r_l\rangle$ and $T_2=\langle y_1,\ldots ,y_p\mid s_1,\ldots ,s_q\rangle$ be presentations. Suppose we have a map $\phi_0\colon F( x_1,\ldots ,x_k)\to F(y_1,\ldots ,y_p)$ that induces a group homomorphism $G(T_1)\to G(T_2)$. For $n>k$ we can extend $\phi_0$ to $\phi\colon F(x_1,\ldots ,x_n)\to F(y_1,\ldots ,y_p, x_{k+1},\ldots ,x_n)$ by defining $\phi(x_i)=x_i$ for $i>k$. Now let $P_1$ be a presentation of the form 
$$P_1=\langle x_1,\ldots ,x_k, x_{k+1},\ldots ,x_n\mid r_1,\ldots ,r_l, r_{l+1},\ldots ,r_m\rangle.$$ We assume each relator $r_j$, $j>l$, contains a generator $x_i$, $i>k$, so that $P_{Y_1}=T_1$ for $Y_1=\{ x_1,\ldots ,x_k\}$. Let
$$P_2=\langle y_1,\ldots ,y_p, x_{k+1},\ldots ,x_n \mid s_1,\ldots ,s_q, \phi(r_{l+1}),\ldots ,\phi(r_m)\rangle.$$
We use $\phi\colon P_1\to P_2$ as shorthand for the situation just described.

\lemanf\label{sbarP} Suppose we have $\phi\colon P_1 \to P_2$. If $P_2$ is DR directed away from $Y_2=\{ y_1,\ldots ,y_p\}$, then $P_1$ is DR directed away from $Y_1=\{ x_1,\ldots ,x_k\}$.
\lemende

See \cite{HaRo19} for a proof. The simplest choice for $T_2$ is the empty presentation $T_2=\langle \ |\ \rangle$, that is $G(T_2)$ is trivial as in the following example.

\bspanf \em Consider $$P_1=\langle x_1,\ldots ,x_k, a, b\mid u_1au_2bu_3a^{-1}u_4b^{-1} \rangle \to P_2= \langle a, b\mid aba^{-1}b^{-1}\rangle,$$
where the $u_i$ are words in $x_1^{\pm 1},\ldots ,x_k^{\pm 1}$. Since $P_2$ is $DR$ (directed away from $\emptyset$), it follows that $P_1$ is DR directed away from $Y=\{ x_1,\ldots ,x_k \}$. More general: Take any DR presentation $P_2$. Add generators $Y=\{ x_1,\ldots ,x_k \}$ and insert words in $Y^{\pm 1}$ into the relators of $P_2$ and one obtains a presentation $P_1$ which is DR directed away from $Y$.
\bspende

\begin{figure}[ht]\centering
\begin{tikzpicture}[scale=1.2]
\fill[red] (0,0) circle (1.6pt); \node[red, below] at (0,0) {$x_1$};
\fill[red] (1,0) circle (1.6pt); \node[red, below] at (1,0) {$x_2$};
\fill[red] (2,0) circle (1.6pt); \node[red, below] at (2,0) {$x_3$};
\fill (3,0) circle (1.6pt); \node[below] at (3,0) {$x_4$};
\fill (4,0) circle (1.6pt); \node[below] at (4,0) {$x_5$};

\begin{scope}[decoration={markings, mark=at position 0.5 with {\arrow{<}}}]
\draw [postaction={decorate}] (3,0) -- (4,0) node[midway, below]{$x_1$};
\end{scope}

\begin{scope}[decoration={markings, mark=at position 0.5 with {\arrow{>}}}]
\draw [red,postaction={decorate}] (0,0) -- (1,0) node[midway, below]{$x_3$};
\draw [red,postaction={decorate}] (1,0) -- (2,0) node[midway, below]{$x_1$};
\draw [postaction={decorate}] (2,0) -- (3,0) node[midway, below]{$x_5$};
\end{scope}

\end{tikzpicture}
\caption{\label{aLOTH} An aspherical  LOT $\Gamma$ with sub-LOT $\Gamma_0$.}
\end{figure}
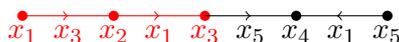

Theorem \ref{sforest} is a strong tool for showing asphericity of LOTs. Here is a simple example:

\bspanf \em Let $\Gamma$ be the LOT depicted in Figure \ref{aLOTH}. Let $\Gamma_0$ be the sub-LOT consisting of the first two edges between $x_1$ and $x_3$ (colored red). $K(\Gamma_0)$ is VA by Theorem \ref{sinj} because $\Gamma_0$ is injective. Since $lk^-(K(\Gamma))$ is a forest relative to $lk^-(K(\Gamma_0))$ Theorem \ref{sforest} implies that $K(\Gamma)$ is VA relative to $K(\Gamma_0)$. By Theorem \ref{sVAsub} $K(\Gamma)$ is VA.
\bspende

Building an example out of an aspherical sub-LOT can even be done in a 2-step (or more steps) process such that the edges you add lead to cycles in the negative and the positive graph:

\bspanf\label{ex2step} \em Let $T=\langle X\mid R\, \rangle$ be an arbitrary compressed VA LOT presentation. Assume $x_i,x_j,x_p,x_q\in X$ (not necessarily pairwise distinct).\\
Let $P_1=\langle X,a,b\mid R, ab=bx_i, ba=ax_j\rangle$. Since $lk^-(K(P_1))$ is a forest relative to $lk^-(K( T))$ Theorem \ref{sforest} implies that $K(P_1)$ is VA relative $K(T)$. Now let
$$P_2=\langle X,a,b,c,d\mid R, ab=bx_i, ba=ax_j, x_pc=cd, x_qd=dc\rangle$$

\begin{figure}[ht]\centering
\begin{tikzpicture}[scale=1.2]
\draw[red] (0,0) circle (1cm); \node[red, left] at (0.2,0) {\Large $T$};
\fill[red] (1,0) circle (1.6pt); \node[red, left] at (1,0) {$x_i$};
\fill[red] (-1,0) circle (1.6pt); \node[red, right] at (-1,0) {$x_j$};
\fill[red] (0,1) circle (1.6pt); \node[red, below] at (0,1) {$x_q$};
\fill[red] (0,-1) circle (1.6pt); \node[red, above] at (0,-1) {$x_p$};

\begin{scope}[decoration={markings, mark=at position 0.5 with {\arrow{>}}}]
\draw [postaction={decorate}] (0,1) -- (0,2) node[midway, right]{$d$};
\draw [postaction={decorate}] (2,0) -- (1,0) node[midway, below]{$b$};
\draw [postaction={decorate}] (0,-1) -- (0,-2) node[midway, right]{$c$};
\draw [postaction={decorate}] (-2,0) -- (-1,0) node[midway, below]{$a$};
\end{scope}

\fill (2,0) circle (1.6pt); \node[above] at (2,0) {$a$};
\fill (-2,0) circle (1.6pt); \node[above] at (-2,0) {$b$};
\fill (0,2) circle (1.6pt); \node[right] at (0,2) {$c$};
\fill (0,-2) circle (1.6pt); \node[right] at (0,-2) {$d$};

\end{tikzpicture}
\caption{\label{aex2step} An aspherical  LOT build in two steps.}
\end{figure}
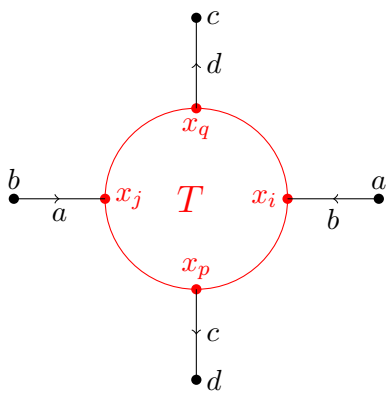
(see Figure \ref{aex2step}).
Since $lk^+(K(P_2))$ is a forest relative to $lk^+(K(P_1))$, Theorem \ref{sforest} and Theorem \ref{sVAsub} imply that $K(P_2)$ is a VA LOT-complex.
\bspende

In this example edges were added to the LOT corresponding to $T$ which have cycles in both, negative and positive graph but still asphericity may be shown by using Theorem \ref{sforest} twice.

\bspanf \em Let $T=\langle X\mid R\, \rangle$ be an arbitrary compressed VA LOT presentation. Assume $x_i,x_j\in X$ (not necessarily distinct). Let 
\[ P=\langle X,u,v,w,y\mid R,  wy=yu, x_iw=wu,  vw=wy, vu=ux_j \rangle\]
(see Figure \ref{aexBP}). Since $lk^+(K(P))$ is a forest relative to $lk^+(K(T))$ Theorem \ref{sforest} implies that $K(P)$ is a VA LOT-complex.

\begin{figure}[ht]\centering
\begin{tikzpicture}[scale=1.2]
\draw[red] (0,0) circle (1cm); \node[red, left] at (0.2,0) {\Large $T$};
\fill[red] (1,0) circle (1.6pt); \node[red, left] at (1,0) {$x_i$};
\fill[red] (-1,0) circle (1.6pt); \node[red, right] at (-1,0) {$x_j$};

\begin{scope}[decoration={markings, mark=at position 0.5 with {\arrow{>}}}]
\draw [postaction={decorate}] (-2,0) -- (-1,0) node[midway, below]{$u$};
\draw [postaction={decorate}] (-3,0) -- (-2,0) node[midway, below]{$w$};
\draw [postaction={decorate}] (1,0) -- (2,0) node[midway, below]{$w$};
\draw [postaction={decorate}] (3,0) -- (2,0) node[midway, below]{$y$};
\end{scope}

\fill (2,0) circle (1.6pt); \node[above] at (2,0) {$u$};
\fill (-2,0) circle (1.6pt); \node[above] at (-2,0) {$v$};
\fill (3,0) circle (1.6pt); \node[above] at (3,0) {$w$};
\fill (-3,0) circle (1.6pt); \node[above] at (-3,0) {$y$};

\end{tikzpicture}
\caption{\label{aexBP} An aspherical  LOT.}
\end{figure}

Let $H$ be the group defined by $T$ and let
\[\hat P=\langle H,u,v,w,y\mid wy=yu, x_iw=wu,  vw=wy, vu=ux_j \rangle\]
be a relative presentation in the sense of Bogley and Pride (see \cite{BP92}). In their paper they define a weight test. $\hat P$ does not satisfy this weight test if $x_i^3=x_j^3$ (which is of course satisfied for $x_i=x_j$). This can be seen by drawing the Whitehead graph in the sense of Bogley and Pride and using the simplex method to show that the weight test does not apply.
\bspende

\bspanf\label{bFWT} \em Let $T=\langle X \mid R\, \rangle$ be an arbitrary compressed aspherical LOT presentation. Assume $x_i,x_j\in X$ (not necessarily distinct). Let 
$$P=\langle X,u,v,w,y,\mid R, uv=vw, x_iy=yw, vw=wy, vu=ux_j \rangle$$ 
(see Figure \ref{aexBP2}).

\begin{figure}[ht]\centering
\begin{tikzpicture}[scale=1.2]
\draw[red] (0,0) circle (1cm); \node[red, left] at (0.2,0) {\Large $T$};
\fill[red] (1,0) circle (1.6pt); \node[red, left] at (1,0) {$x_i$};
\fill[red] (-1,0) circle (1.6pt); \node[red, right] at (-1,0) {$x_j$};

\begin{scope}[decoration={markings, mark=at position 0.5 with {\arrow{>}}}]
\draw [postaction={decorate}] (-2,0) -- (-1,0) node[midway, below]{$u$};
\draw [postaction={decorate}] (-2,0) -- (-3,0) node[midway, below]{$w$};
\draw [postaction={decorate}] (1,0) -- (2,0) node[midway, below]{$y$};
\draw [postaction={decorate}] (3,0) -- (2,0) node[midway, below]{$v$};
\end{scope}

\fill (2,0) circle (1.6pt); \node[above] at (2,0) {$w$};
\fill (-2,0) circle (1.6pt); \node[above] at (-2,0) {$v$};
\fill (3,0) circle (1.6pt); \node[above] at (3,0) {$u$};
\fill (-3,0) circle (1.6pt); \node[above] at (-3,0) {$y$};

\end{tikzpicture}
\caption{\label{aexBP2} An aspherical  LOT.}
\end{figure}
 $P$ and $T$ do not satisfy the conditions of Theorem \ref{sforest} since cycles occur in $lk^+(K(P))$ relative to $lk^+(K(T))$ and $lk^-(K(P))$ relative to $lk^-(K(T))$. There is also no way to build up $P$ in two steps as in Example \ref{ex2step}. On the other hand it can easily be seen that $P$ satisfies the weight test of Bogley and Pride (see \cite{BP92}). Give all edges which occur twice in $lk(K(P-T))$ weight 1 and all other edges weight 0. 
$P$ satisfies the relative weight test of Theorem \ref{sWTb} since it is injective relative to $T=\langle {X}\mid R \rangle$. See page 11, after Lemma 5 in \cite{HaRo17} for the definition of ``injective relative to''.
\bspende

The small cancellation conditions C($p$), T($q$) are defined for instance in the book of Lyndon and Schupp \cite{LS77}. The following theorem is shown in \cite{HaRo19}. It is an application of the weight test (Theorem \ref{sWT}):

\satzanf\label{s44} Let $P$ be a finite presentation with cyclically reduced relators and $Y$ a subset of the generators. Assume that $P$ is C(4), T(4) or C(6), T(3) and that no two consecutive letters in a (cyclically read) relator of $P$ are elements of $Y^{\pm 1}$. Then $P$ is DR directed away from $Y$.
\satzende

\begin{figure}[ht]\centering
\begin{tikzpicture}
\draw (1,0)--(2,0)--(3,0)--(4,0)--(5,0)--(6,0)--(7,0);
\fill (1,0) circle (2pt);
\fill (2,0) circle (2pt);
\fill (3,0) circle (2pt);
\fill (4,0) circle (2pt);
\fill (5,0) circle (2pt);
\fill (6,0) circle (2pt);
\fill (7,0) circle (2pt);
\node [below] at (1,0) {$x_1$};
\node [above] at (1.5,0) {$x_3$};
\node [below] at (2,0) {$x_2$};
\node [above] at (2.5,0) {$x_4$};
\node [below] at (3,0) {$x_3$};
\node [above] at (3.5,0) {$x_5$};
\node [below] at (4,0) {$x_4$};
\node [above] at (4.5,0) {$x_6$};
\node [below] at (5,0) {$x_5$};
\node [above] at (5.5,0) {$x_7$};
\node [below] at (6,0) {$x_6$};
\node [above] at (6.5,0) {$x_1$};
\node [below] at (7,0) {$x_7$};
\end{tikzpicture}
\caption{\label{alc4t4} A labelled oriented tree which is C(4), T(4).}
\end{figure}
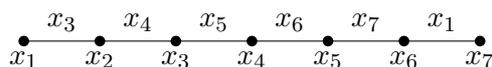

In Rosebrock \cite{Ro94} it is described how to check whether a LOT is C(4), T(4).
If a LOT is C(4), T(4) then Theorem \ref{s44} implies that it is DR away from any of its generators. But there is more in concrete cases:

\bspanf \em Consider the LOT $\Gamma$ of Figure \ref{alc4t4} with any orientation of its edges. This LOT is C(4), T(4). If you choose $Y$ to be the edges of $K(\Gamma)$ corresponding to one of the sets
\[\{ x_1,x_2,x_5\} , \{ x_1,x_2\} , \{ x_1,x_4\} , \{ x_1,x_5\} , \{ x_2,x_5\} , \{ x_2,x_6\} , \{ x_2,x_7\} , \{ x_3,x_6\} ,\]
\[ \{ x_3,x_7\} , \{ x_4,x_7\} \]
then Theorem \ref{s44} implies that $K(\Gamma)$ is DR away from $Y$.
\bspende

Let $\Gamma$ be a compressed LOT and $\Gamma_0$ be a maximal proper sub-LOT with vertex set $Y$. Let $y\in\Gamma_0$ be any vertex. Let $\bar{\Gamma}$ be the LOT obtained from $\Gamma$ by collapsing all of $\Gamma_0$ to the vertex $y$. Every occurrence of a vertex $x$ of $\Gamma_0$ in $\Gamma -\Gamma_0$ is replaced by $y$ in $\bar{\Gamma}$.  The following theorem is shown in \cite{HaRo19}, it is an application of Lemma \ref{sbarP}, Theorem \ref{sminmax} and Theorem \ref{s44}:

\satzanf\label{sLOT} If $\bar\Gamma$ is compressed and $K(\bar\Gamma )$ is DR directed away from $y$ then $\Gamma$ is DR directed away from $Y$, the edges of $\Gamma_0$. In particular if
\begin{enumerate}
\item $lk^+(K(\bar \Gamma))$, or $lk^-(K(\bar \Gamma))$ is a tree, or
\item $lk(K(\bar \Gamma))$ does not contain cycles of length less than four,
\end{enumerate}
then $\Gamma$ is DR directed away from $Y$.
\satzende

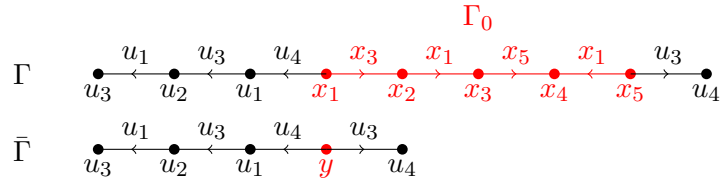
\begin{figure}[ht]\centering
\begin{tikzpicture}
\node at (-1,0) {$\Gamma$};
\node[red] at (5,0.75) {$\Gamma_0$};

\fill (0,0) circle (2pt);
\node[below] at (0,0) {$u_3$};
\fill (1,0) circle (2pt);
\node[below] at (1,0) {$u_2$};
\fill (2,0) circle (2pt);
\node[below] at (2,0) {$u_1$};
\fill[red] (3,0) circle (2pt);
\node[red,below] at (3,0) {$x_1$};
\fill[red] (4,0) circle (2pt);
\node[red,below] at (4,0) {$x_2$};
\fill[red] (5,0) circle (2pt);
\node[red,below] at (5,0) {$x_3$};
\fill[red] (6,0) circle (2pt);
\node[red,below] at (6,0) {$x_4$};
\fill[red] (7,0) circle (2pt);
\node[red,below] at (7,0) {$x_5$};
\fill (8,0) circle (2pt);
\node[below] at (8,0) {$u_4$};

\begin{scope}[decoration={markings, mark=at position 0.5 with {\arrow{>}}}]
\draw [red,postaction={decorate}] (3,0) -- (4,0) node[midway, above]{$x_3$};
\draw [red,postaction={decorate}] (4,0) -- (5,0) node[midway, above]{$x_1$};
\draw [red,postaction={decorate}] (5,0) -- (6,0) node[midway, above]{$x_5$};
\draw [postaction={decorate}] (7,0) -- (8,0) node[midway, above]{$u_3$};
\end{scope}

\begin{scope}[decoration={markings, mark=at position 0.5 with {\arrow{<}}}]
\draw [postaction={decorate}] (0,0) -- (1,0) node[midway, above]{$u_1$};
\draw [postaction={decorate}] (1,0) -- (2,0) node[midway, above]{$u_3$};
\draw [postaction={decorate}] (2,0) -- (3,0) node[midway, above]{$u_4$};
\draw [red,postaction={decorate}] (6,0) -- (7,0) node[midway, above]{$x_1$};
\end{scope}

\begin{scope}[yshift=-1cm]
\node at (-1,0) {$\bar{\Gamma}$};

\fill (0,0) circle (2pt);
\node[below] at (0,0) {$u_3$};
\fill (1,0) circle (2pt);
\node[below] at (1,0) {$u_2$};
\fill (2,0) circle (2pt);
\node[below] at (2,0) {$u_1$};
\fill[red] (3,0) circle (2pt);
\node[red,below] at (3,0) {$y$};
\fill (4,0) circle (2pt);
\node[below] at (4,0) {$u_4$};

\begin{scope}[decoration={markings, mark=at position 0.5 with {\arrow{>}}}]
\draw [postaction={decorate}] (3,0) -- (4,0) node[midway, above]{$u_3$};
\end{scope}

\begin{scope}[decoration={markings, mark=at position 0.5 with {\arrow{<}}}]
\draw [postaction={decorate}] (0,0) -- (1,0) node[midway, above]{$u_1$};
\draw [postaction={decorate}] (1,0) -- (2,0) node[midway, above]{$u_3$};
\draw [postaction={decorate}] (2,0) -- (3,0) node[midway, above]{$u_4$};
\end{scope}

\end{scope}

\end{tikzpicture}
\caption{\label{aex1} The LOT $\Gamma$ with sub-LOT $\Gamma_0$ and the LOT $\bar{\Gamma}$}
\end{figure}

\bspanf \em Figure \ref{aex1} shows a compressed LOT $\Gamma$ together with a sub-LOT $\Gamma_0$ with vertices $Y=\{ x_1,\ldots ,x_5\}$. Below $\Gamma$ we see the LOT $\bar{\Gamma}$ obtained from $\Gamma$ by collapsing $\Gamma_0$ to the vertex $y$.
\bspende

\vspace{-1ex}

\noindent Note that $lk^+(K(\bar{\Gamma}))$ is a tree, so by Theorem \ref{sLOT} $\Gamma$ is DR directed away from $Y$. Also observe that $lk^-(K(\Gamma_0))$ is a tree, so $K(\Gamma_0)$ is DR, which implies that $K(\Gamma )$ is aspherical by Theorem \ref{RelPiS} (1). Note that neither $lk^+(K(\Gamma ))$ nor $lk^-(K(\Gamma ))$ is a tree.\\

The process of collapsing a sub-LOT in a given labeled oriented tree can also be reversed: If $\bar {\Gamma}$ is a labeled oriented tree, $y$ is a vertex in $\bar {\Gamma}$, and $\Gamma_0$ is a labeled oriented tree, we remove $y$ from $\bar {\Gamma}$ and insert $\Gamma_0$ to obtain a labeled oriented tree $\Gamma$ that contains $\Gamma_0$. Collapsing $\Gamma_0$ in $\Gamma$ to a vertex $y$  brings us back to $\bar {\Gamma}$. So the previous theorem can also be stated as follows: If $\bar{\Gamma}$ is a labeled oriented tree that satisfies either condition (1) or (2) of Theorem \ref{sLOT}, then inserting any LOT $\Gamma_0$ into $\bar{\Gamma}$ results in a labeled oriented tree $\Gamma$ for which  $K(\Gamma)$ is DR directed away from the set $Y$ of edges of $K(\Gamma_0)$.\\

LOTs $\Gamma$ where $lk^+(K({\Gamma}))$ or $lk^-(K({\Gamma}))$ is a tree abound. 
If  a LOT $\Gamma'$ is obtained from a LOT $\Gamma$ by changing some edge orientations, we call $\Gamma'$ a {\em reorientation} of $\Gamma$. In \cite{HR01} Proposition 5.1, Huck and Rosebrock show that each LOT $\Gamma$ has a reorientation $\Gamma'$ such that $lk^+(K(\Gamma'))$ is a tree. Theorem \ref{sminmax} now implies

\satzanf
Each LOT $\Gamma$ has a reorientation $\Gamma'$ so that $K(\Gamma')$ is DR away from any one of its edges.
\satzende

In \cite{Ro94} the second author gives conditions on a labeled oriented tree so that condition (2) of Theorem \ref{sLOT} holds. 

\begin{figure}[ht]\centering
\begin{tikzpicture}
\node at (-1,0) {$\Gamma$};
\node[red] at (5,0.75) {$\Gamma_0$};

\fill (0,0) circle (2pt);
\node[below] at (0,0) {$u_1$};
\fill (1,0) circle (2pt);
\node[below] at (1,0) {$u_2$};
\fill (2,0) circle (2pt);
\node[below] at (2,0) {$u_3$};
\fill[red] (3,0) circle (2pt);
\node[red,below] at (3,0) {$u_4$};
\fill(7,0)[red] circle (2pt);
\node[red, below] at (7,0) {$u'_4$};
\fill (8,0) circle (2pt);
\node[below] at (8,0) {$u_5$};
\fill (9,0) circle (2pt);
\node[below] at (9,0) {$u_6$};
\fill (10,0) circle (2pt);
\node[below] at (10,0) {$u_7$};

\draw [red,postaction={decorate}] (3,0) -- (4,0); 
\draw [red,postaction={decorate}] (4,0) -- (5,0);
\draw [red,postaction={decorate}] (5,0) -- (6,0);
\draw [postaction={decorate}] (7,0) -- (8,0);
\draw [postaction={decorate}] (8,0) -- (9,0);
\draw [postaction={decorate}] (9,0) -- (10,0);

\draw [postaction={decorate}] (0,0) -- (1,0) node[midway, above]{$u_3$};
\draw [postaction={decorate}] (1,0) -- (2,0) node[midway, above]{$u_4$};
\draw [postaction={decorate}] (2,0) -- (3,0) node[midway, above]{$u_5$};
\draw [red,postaction={decorate}] (6,0) -- (7,0);
\draw [postaction={decorate}] (7,0) -- (8,0) node[midway, above]{$u_6$};
\draw [postaction={decorate}] (8,0) -- (9,0) node[midway, above]{$u_7$};
\draw [postaction={decorate}] (9,0) -- (10,0) node[midway, above]{$u_1$};

\end{tikzpicture}
\caption{\label{aex2} A labeled oriented tree $\Gamma$ with sub-LOT $\Gamma_0$.}
\end{figure}
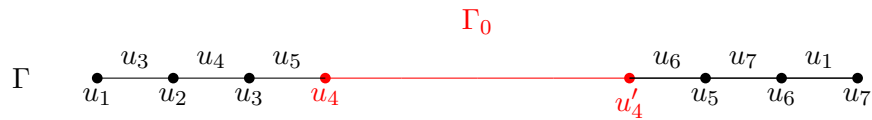

\bspanf \em Figure \ref{aex2} shows a labelled oriented tree $\Gamma$ (orientations can be chosen at will) with a sub-LOT $\Gamma_0$ between $u_4$ and $u_4'$ (which can be filled in at will).
\bspende

\vspace{-1ex}

\noindent Notice that if we collapse the red sub-LOT $\Gamma_0$ to the vertex $y=u_4$, we obtain a labeled oriented tree $\bar{\Gamma}$ for which $lk(K(\bar{\Gamma}))$ does not contain cycles of length less than four. It follows from Theorem \ref{sLOT} that $K(\Gamma)$ is DR away from the edge set of $K(\Gamma_0)$.\\

At last we provide an example to Theorem \ref{sWTa}.

\bspanf  \em Let $n\ge 2$ and $1\le k\le n$. Let $P_1=\langle y_1,\ldots ,y_k\mid R_1\,\rangle$ and
$P_2=\langle y_{k+1},\ldots ,y_n\mid R_2\,\rangle$ and $K_i$ be the standard 2-complex corresponding to $P_i$.  Assume  $R_i$ is a set of relators, such that each $y_m$ is nontrivial in $\pi_1(K_i)$ $(i=1,2)$. Let $L$ be the standard 2-complex given by
\[ P=\langle x_1,\ldots ,x_n,y_1,\ldots ,y_n\mid R_1,R_2, w_1,\ldots ,w_n\rangle\]
where $w_i=x_i^{-2}y_ix_{i+1}y_i^{-1}$ $(i \mbox{ mod }n)$.
By Theorem \ref{sWTa} $L$ is VA relative $K=K_1\vee K_2$. This can be seen by assigning weight 0 to corners in $lk^+(L-K)$ or $lk^-(L-K)$ and weight 1 for the other corners of $L-K$. Observe that relators have not exponent sum 0.
\bspende

\vspace{3ex}

\noindent Jens Harlander\\
Department of Mathematics\\
Boise State University\\
Boise, ID 83725-1555\\
USA\\
email: jensharlander@boisestate.edu\\

\noindent Stephan Rosebrock\\
P{\"a}dagogische Hochschule Karlsruhe\\
Bismarckstr. 10\\
76133 Karlsruhe\\
Germany\\
email: rosebrock@ph-karlsruhe.de


\begin{thebibliography}{99}

\bibitem{BM19}
M. A. Blufstein and E. G. Minian,
\newblock{\it Strictly systolic angled complexes and hyperbolicity of one-relator groups}, 
https://arxiv.org/abs/1907.06738, (2019).

\bibitem{BP92}
W.~Bogley and S.~Pride,
\newblock {\it Aspherical relative presentations.}
\newblock Proc.~of the Edinburgh Math.~Society 35, (1992), pp. 1--39.

\bibitem{BEW18}
W.~Bogley, M.~Edjvet and G.~Williams,
\newblock {\it Aspherical relative presentations all over again}, in: Groups St Andrews 2017 in Birmingham, London Mathematical Society Lecture Note Series, Editors: Campbell, C. M. and Parker, C. W. and Quick, M. R. and Robertson, E. F. and Roney-Dougal, Cambridge University Press (2019), pp. 169-199.

\bibitem{CT96} J. M. Corson and B. Trace, {\it Geometry and algebra of nonspherical 2-complexes}, J. London Math. Soc. 54, (1996), pp. 180-198.

\bibitem{CT00}
J. M. Corson and B. Trace, {\it Diagrammatically reducible complexes and Haken manifolds},
J. Austral. Math. Soc. (Series A) 69, (2000), pp. 116-126.

\bibitem{Ger87} S. M. Gersten, {\it Reducible diagrams and equations over groups},
in: Essays in Group Theory (S. M. Gersten editor), MSRI Publications 8 (1987), pp.~15--73.

\bibitem{HaRo17} J. Harlander and S. Rosebrock, {\em Injective labelled oriented trees are aspherical}, Mathematische Zeitschrift 287 (1), (2017), pp. 199--214.

\bibitem{HaRo19} J. Harlander and S. Rosebrock, {\em Directed diagrammatic reducibility}, to appear in Topology and its Applications (2020).

\bibitem{HaRo19b} J. Harlander and S. Rosebrock, {\em Relative Vertex Asphericity}, to appear in: Canad. Math. Bull. (2020).

\bibitem{HR01}
G. Huck and S. Rosebrock, {\em Aspherical Labelled Oriented Trees and Knots}, Proceedings of the Edinburgh Math.~Soc. 44, (2001), pp. 285-294.

\bibitem{LS77}
R.~Lyndon and P.~Schupp.
\newblock {\it Combinatorial group theory}, Springer Verlag, Berlin, (1977).

\bibitem{Ro18} S. Rosebrock, {\em Labelled Oriented Trees and the Whitehead Conjecture}; 
Advances in Two-Dimensional Homotopy and Combinatorial Group Theory;
Cambridge University Press, LMS Lect.~Notes  446, editors
W.~Metzler, S.~Rosebrock; (2018), pp. 72--97.

\bibitem{Ro94} S. Rosebrock, {\em On the realization of Wirtinger presentations as knot groups}, Journal of Knot Theory and its Ramifications, Vol. 3 (1994), pp. 211-222.

\bibitem{S83} A. Sieradski, {\em A coloring test for asphericity}, Quart. J. Math. Oxford (2) 34 (1983), pp. 97--106.

\end{thebibliography}
\end{document}